\documentclass[final]{siamltex}

\usepackage{pdfsync}
\usepackage{graphicx}
\usepackage{amssymb}
\usepackage{amsfonts}
\usepackage{mathrsfs}
\usepackage{marvosym}
\usepackage[pdftex]{color}

\usepackage[dvips]{epsfig} \DeclareGraphicsExtensions{.ps,.eps}
\usepackage{url}

\newtheorem{remark}[theorem]{Remark}

\def\Z{\mathbb{Z}}
\def\R{\mathbb{R}}
\def\C{\mathbb{C}}
\def\T{\mathbb{T}}
\def\D{\mathbb{D}}

\def\H{\mathbb{H}}
\def\P{\mathbb{P}}
\newcommand{\dr}[2]{\D_r(#1\|#2)}
\def\im{{\rm im}\,}
\def\det{{\rm det}\,}
\def\exp{{\rm exp}\,}

\def\tr{{\rm tr}\,}
\newcommand{\nn}{\nonumber}
\newcommand{\imunit}{{\rm j}}
\newcommand{\tp}{^{\top}}

\newcommand{\eq}{\begin{equation}}
\newcommand{\eeq}{\end{equation}}
\newcommand{\eqn}{\begin{eqnarray}}
\newcommand{\eeqn}{\end{eqnarray}}
\newcommand{\Hermitian}{\H}
\newcommand{\Cthm}{C( \mathbb{T},\Hermitian_m)}

\newcommand{\pjb}{\Pi_B}

\newcommand{\Rgamma}{\mathop{\rm Range} \Gamma}

\newcommand{\Tr} {\mbox{\rm tr}}

\newcommand{\qed}{\hfill $\Box$ \vskip 2ex}

\def\bmat{\left[ \begin{array}}
\def\emat{\end{array} \right]}

\newcommand{\bsea}{\begin{subeqnarray}}
\newcommand{\esea}{\end{subeqnarray}}

\def\bmat{\left[ \begin{array}}
\def\emat{\end{array} \right]}

\newcommand{\h}[1]{h_r({#1})}
\newcounter{acount}

\definecolor{Royalblue}{cmyk}{1,0.30,0.2,0.2}

\newfont{\BB}{msbm10 scaled\magstep1}
\def\E{\mbox{\BB E}}

\title{On the Geometry of Maximum Entropy Problems\thanks{Work partially supported by the Italian Ministry for Education and Resarch (MIUR) under PRIN grant
``Identification and Robust Control of Industrial Systems", by the CPDA080209/08 and QFUTURE research grants of the University of Padova and by the Department of Information Engineering research project ``QUINTET".}}

\author{Michele Pavon \thanks{Dipartimento di Matematica, Universit\`a di Padova, via
Trieste 63, 35131 Padova, Italy ({\tt
pavon@math.unipd.it}).}
        \and Augusto Ferrante \thanks{Dipartimento di Ingegneria dell'Informazione, Universit\`a di Padova,
via Gradenigo 6/B, 35131
Padova, Italy. ({\tt  augusto@dei.unipd.it})}}

\begin{document}

\maketitle

\begin{abstract}
We show that  a  simple geometric result  suffices to derive the form of the optimal solution in a large class of finite and infinite-dimensional {\em maximum entropy problems} concerning probability distributions, spectral densities and covariance matrices. These include {\em Burg's spectral estimation method} and  {\em Dempster's covariance completion}, as well as various recent generalizations of the above. We then apply this orthogonality principle to the new problem of completing a block-circulant covariance matrix when an {\em a priori} estimate is available.
\end{abstract}



\pagestyle{myheadings}
\thispagestyle{plain}
\markboth{MICHELE PAVON AND AUGUSTO FERRANTE}{GEOMETRY OF MAXIMUM ENTROPY PROBLEMS}

\section{Prelude: Four famous maximum entropy problems}
In this section, we briefly review four classical maximum entropy problems that have played an important role in the history of various scientific areas. These are namely problems where entropy is maximized under linear constraints. We shall later derive the form of the optimal solution in three of these problems by the same geometric principle (Theorem \ref{hilbertspace} in Section \ref{geometricresult}).

\subsection{1877: Boltzmann's loaded dice} In 1877, Boltzmann \cite[p.169]{BOL} posed the following question: Consider $N$ molecules that can only take the following $p+1$ values of kinetic energy \footnote{``lebendige Kraft", the classical {\em vis viva} originating with Gottfried Leibniz which was actually twice the kinetic energy.} $0,\epsilon,2\epsilon,\ldots,p\epsilon$. Suppose $n_i$ molecules have kinetic energy $i\epsilon, i=0,1,\ldots,p$. We then have a ``macrostate", a ``Zustandverteilung" in Boltzmann's language\footnote{the expression ``Komplexion" in \cite{BOL} refers instead to a {\em microstate} and not to a macrostate as stated in \cite[Section 4]{uffink}.}, indexed by $(n_0,n_1,\ldots,n_p)$ corresponding to the  multinomial coefficient  
$$\frac{N!}{n_0!n_1!\ldots n_p!}$$ 
``microstates" each having probability $(p+1)^{-N}$. Suppose that the sum of the kinetic energy of all molecules is a given quantity $\lambda\epsilon=L$. Boltzmann proceeded to find the macrostate which corresponds to more microstates, namely that has highest probability,  among those having total kinetic energy $L$. This  is, to the best of our knowledge,  the first maximum entropy problem in history.

Boltzmann's problem was popularized in the following form \cite{Jaynes82,COVER_THOMAS}. Suppose $N$ dice are rolled and we are informed that the total number of spots is $N\cdot 4.5$. We are asked: What proportion of the dice are showing face $i, i=1,2,\ldots,6$? The number of different ways that $N$ dice can fall so that $n_i$ dice show face $i$ is given by 
\begin{equation}\label{multinomialcoeff}
\frac{N!}{n_1!n_2!\ldots n_6!},\quad \sum_{i=1}^6 n_i=N.
\end{equation}
Again, the ``macrostate" $(n_1,n_2,\ldots,n_6)$ corresponds to $\frac{N!}{n_1!n_2!\ldots n_6!}$ ``microstates" each having probability $6^{-N}$. To find the most probable macrostate, we need to maximize the multinomial coefficient (\ref{multinomialcoeff})  under the constraint
\begin{equation}\label{meanconstraint}\sum_{i=1}^6 i\cdot n_i=N\cdot 4.5.
\end{equation}
This procedure will yield the macrostate, among those satisfying (\ref{meanconstraint}), that {\em can be realized in more ways}. Assuming that $N$ is large, we now use a crude version of Stirling's approximation $N!\approx e^{-N}N^N$. We get
\begin{eqnarray}\nonumber &&\frac{N!}{n_1!n_2!\ldots n_6!}\approx\frac{e^{-N}N^N}{\prod_{i=1}^6 e^{-n_i} n_i^{n_i}}=\prod_{i=1}^6\left(\frac{N}{n_i}\right)^{n_i}=\prod_{i=1}^6 e^{-n_i\ln\left(\frac{n_i}{N}\right)}=\\&&e^{-\sum_{i=1}^6n_i\ln\left(\frac{n_i}{N}\right)}=e^{NH(p)},\quad p_i=\frac{n_i}{N}, i=1,2,\ldots,6.\nonumber
\end{eqnarray}
Thus, for $N$ large, maximizing (\ref{multinomialcoeff}) under (\ref{meanconstraint}) is almost equivalent to maximizing the entropy
$$H(p)=-\sum_{i=1}^6 p_i\ln\left(p_i\right)
$$
under the  constraint 
\begin{equation}\label{newconstraint}\sum_{i=1}^6 i\cdot p_i=4.5.
\end{equation}
The solution has the form
\begin{equation}\label{BOLT}p_i^*=\frac{e^{\lambda_i}}{\sum_{i=1}^6 e^{\lambda_i}},
\end{equation}
where the $\lambda_i$ must be such that
$$\sum_{i=1}^6 i\cdot \frac{e^{\lambda_i}}{\sum_{i=1}^6 e^{\lambda_i}}=4.5.
$$
Hence, the most probable macrostate is $(Np_1^*,Np_2^*,\ldots,Np_6^*)$ and we expect to find $n_i^*=Np_i^*$ dice showing face $i$. More is true: It can be shown \cite[Chapter 13]{COVER_THOMAS} that, for $N$ large, with probability close to one, other distributions satisfying (\ref{newconstraint}) are close to $p^*$. This fact is sometimes referred to as Entropy Concentration Theorem \cite{Jaynes82}. More generally, when $ F(p):=-\sum_k p_k\log(p_k)$, the maximizer of $F$ subject to a linear constraint $Lp=c$ has the form of a Boltzmann-Gibbs distribution
\begin{equation}\label{gibbs}
p_k=\frac{1}{Z} e^{-\langle \Lambda, L_k\rangle}
\end{equation}
where $L_k$ is the $k$th column of the matrix $L$ and $Z$ a normalizing constant ({\em partition function}). This can of course also be formulated in the continuous setting (with integrals) and is also a basic result in statistics \cite{csiszar2,CsiszarMatus,CsiszarMatus2}. 

\subsection{1931: Schr\"{o}dinger's Bridges}\label{schroedinger}
In 1931/32, before the very foundations of probability were laid, Erwin Schr\"{o}dinger studied the following abstract problem \cite{S,S2}. Consider the evolution of a cloud of $N$ independent Brownian particles. Here $N$ is large, say of the order of Avogadro's number. This cloud of particles has been observed having
at some initial time $t_0$ an empirical distribution equal to $\rho_0(x)dx$. At some later time $t_1$, an empirical distribution equal to  $\rho_1(x)dx$ is observed which considerably differs from what it should be according to the law of large numbers, namely
$$\left(\int_{t_0}^{t_1}p(t_0,y,t_1,x)\rho_0(y)dy\right)dx,
$$
where
$$p(s,y,t,x)=\left[2\pi(t-s)\right]
^{-\frac{n}{2}}\exp\left[-\frac{|x-y|^2} {2(t-s)}\right],\quad s<t
$$
is the transition density of the Wiener process. It is apparent that the particles have been transported in an unlikely way. But of the many unlikely ways in which this could have happened, which one is
the most likely? 
Schr\"{o}dinger showed that the solution, namely the bridge from $\rho_0$ to $\rho_1$ over Brownian motion, has at each time a density $q$ that factors as $q(x,t)=\varphi(x,t)\hat{\varphi}(x,t)$, where $\varphi$ and $\hat{\varphi}$ solve the system
\begin{eqnarray}\label{SY1}
&&\varphi(t,x)=\int
p(t,x,t_1,y)\varphi(t_1,y)dy,\quad \varphi(t_0,x)\hat{\varphi}(t_0,x)=\rho_0(x)\\&&\hat{\varphi}(t,x)=\int
p(t_0,y,t,x)\hat{\varphi}(t_0,y)dy,\quad \varphi(t_1,x)\hat{\varphi}(t_1,x)=\rho_1(x).\label{SY2}
\end{eqnarray}
It took more than fifty years before F\"{o}llmer, recovering Schr\"{o}dinger's original motivation, observed in \cite{F2} that this is a problem of {\em large deviations\footnote{Large deviations theory has various  applications in hypothesis testing, rate distortion theory, etc, see e.g. \cite[Chapter 11]{COVER_THOMAS}, \cite{DS}, \cite[Chapters 2,3,7]{DZ}. For large deviations of the empirical distribution (level-2 large deviations) for diffusion processes see \cite{F2,Wak} (see also \cite{PT} for a recent extension of this theory to discrete-time classical and quantum evolutions).} of the empirical distribution} on path space \cite{ellis} connected, thanks to Sanov's theorem \cite{SANOV}, to a maximum entropy problem. Schr\"{o}dinger's problem may be considerably generalized. Let $\Omega:={\cal
C}([t_0,t_1],\R^n)$ denote the family of $n$-dimensional
continuous functions, let
$W_x$ denote Wiener measure on $\Omega$ starting at $x$, and let
$$W:=\int_{\R^n} W_x\,dx
$$
be stationary Wiener measure. Let ${\cal D}$ be the family of
distributions on $\Omega$ that are equivalent to $W$.  For
$Q,P\in{\cal D}$, we define  the {\it relative
entropy } $\D(P\|Q)$ of $P$ with respect to $Q$ as
$$\D(P\|Q)=E_P\left[\log\frac{dP}{dQ}\right],$$
where $dP/dQ$ is the Radon-Nikodym derivative of $P$ with respect to $Q$. Let ${\cal D}(\rho_0,\rho_1)$ be distributions in ${\cal D}$ having the observed densities at times $t_0$ and $t_1$.  If there is
at least one
$P$ in
${\cal D}(\rho_0,\rho_1)$ such that
$\D(P\|Q)<\infty$,
it may be shown that there exists a unique minimizer $P_c$ in
${\cal D}(\rho_0,\rho_1)$ called in the language of Csisz\'ar the I-{\em projection} of $Q$ onto ${\cal D}(\rho_0,\rho_1)$ \cite{csiszar0,csiszar1,CsiszarMatus2}. It is  {\em the Schr\"{o}dinger bridge} from $\rho_0$ to $\rho_1$ over $Q$.  In \cite{DGW}, using a conditional version of Sanov's theorem established by Csisz\'ar \cite{csiszar1}, it was shown that such I- projection $P_c$ provides the answer to Schr\"{o}dinger's original question: Namely, the asymptotic empirical distribution on path space, conditioned that the initial and final empirical distributions are $\rho_0(x)dx$ and $\rho_1(y)dy$, respectively, is indeed given by $P_c$.

\subsection{1967: Burg's spectral estimation method}
Suppose the {\em covariance lags} $c_k=\E[y(k)y(0)], k=0,1,\ldots,n-1$ of a stationary, zero-mean, Gaussian process have been estimated from the data. How should one extend the covariance? In 1967, while working on spectral estimation for geophysical data \cite{BURG1}, Burg suggested the following approach. Rather than setting the other covariance lags to zero, one should set them to values such that they maximize the {\em entropy rate} (see Section \ref{matricial functions} below) of the process. The solution is an {\em autoregressive process} of the form
$$y(m)=\sum_{k=1}^{n-1}a_ky(m-k)+w(m),
$$
where $w$ is a zero-mean, Gaussian white noise sequence with variance $\sigma^2$. The parameters $a_1,\ldots,a_{n-1},\sigma^2$ are such that the first $n$ covariance lags match the given ones.

\subsection{1972: Dempster's covariance selection}\label{dempster}
In the seminal paper \cite{Dempster-72},  a general strategy for completing a partially specified covariance matrix was introduced. 
Consider a zero-mean, multivariate Gaussian distribution with density
$$p(x)=(2\pi)^{-n/2}|\Sigma|^{-1/2}\exp\left\{-\frac{1}{2}x\tp\Sigma^{-1}x\right\},\quad x\in\R^n.
$$
Suppose that the elements $\{\sigma_{ij};1\le i\le j\le n, (i,j)\in \bar{{\cal I}}\}$, with $(i,i)\in\bar{{\cal I}}$ for all $i=1\dots n$, have been specified. How should $\Sigma$ be completed? Dempster resorts to a form of the {\em Principle of Parsimony}  of parametric model fitting: As the elements $\sigma^{ij}$ of $\Sigma^{-1}$ appear as natural parameters of the model, one should set  $\sigma^{ij}$ to zero for $1\le i\le j\le n, (i,j)\not\in \bar{{\cal I}}$. Notice that $\sigma^{ij}=0$ has the probabilistic interpretation  that the $i$-th and $j$-th components of the Gaussian random vector  are {\em conditionally independent} given the other components \cite{SK}. 
We say that a positive definite completion $\Sigma^\circ$ of $\Sigma$ is a {\em Dempster Completion} if $[(\Sigma^\circ)^{-1}]_{i,j}=0$ for all $(i,j)\not\in \bar{{\cal I}}$.
In particular, Dempster  proved that when a symmetric, positive-definite completion of $\Sigma$ exists, then there exists a unique Dempster's Completion $\Sigma^\circ$. This completion maximizes the (differential) entropy  
\begin{equation}\label{entropy}
H(p)=-\int_{\R^n}\log (p(x))p(x)dx=\frac{1}{2}\log(\det\Sigma)+\frac{1}{2}n\left(1+\log(2\pi)\right)
\end{equation}
among zero-mean Gaussian distributions having the prescribed elements $\{\sigma_{ij};1\le i\le j\le n, (i,j)\in \bar{{\cal I}}\}$.
Thus, Dempster's Completion $\Sigma^\circ$ solves a {\em maximum entropy} problem, i.e. maximizes entropy under linear constraints. 

\section{Overture}
The long tale of maximum entropy problems originates more than one hundred and thirty years ago with Boltzmann \cite{BOL} at the dawn of statistical mechanics. Since then, several deep thinkers such as Jaynes \cite{Jaynes57,Jaynes82}, Dempster  \cite{Dempster-72},  Csisz\'{a}r \cite{csiszar2}, to name but a few, have tried to  explain the rationale behind the maximum entropy approach. Yet, this method, although never presented as a panacea \cite[p.939]{Jaynes82}, is  still often viewed as an idiosyncratic choice. Before we get all tangled up with ``predict states that can be realized by Nature in the greatest number of ways, while agreeing with your macroscopic information" (Jaynes interpreting Gibbs), apply the Principle of Parsimony of parametric model fitting (Dempster) or the axiomatic approach (Csisz\'ar), we hasten reassure the reader: We are not going to give here even a pr\'ecis of the motivation behind maximum entropy problems. Others have done it much better than we ever could. The scope of this paper is much more modest and yet, in a way, ambitious. 

We want to point out that behind an endless string of maximum entropy solutions there is a simple geometric principle. Namely, that a whole class of seemingly unrelated results concerning probability distributions, spectral densities and covariance matrices are consequences of the {\em same variational principle}.
All these problems feature linear constraints which determine an {\em affine} subspace ${\cal W}$ in which the solution must be sought. Theorem \ref{hilbertspace} (or its generalization Theorem \ref{banachspace}) simply states that the  gradient (or a suitable generalization of it) of the entropy functional at a critical point must belong to the orthogonal complement (or, more generally, to the {\em annihilator}) of the subspace ${\cal V}$ of which the affine space ${\cal W}$ is a translation.  Just to avoid any misunderstanding: We are not dealing here with the (usually challenging) {\em existence problem} \cite{BL,BL2,LEO,LEO2}. We simply want to derive in the most economic way the form of the optimal solutions assuming that they exist.

This orthogonality result is actually a direct consequence of a Lagrange multipliers argument. Nevertheless, we show that when the constraints are linear, there is no need to bring in our illustrious compatriot's multipliers be they vectors, matrices or signals. One can simply skip the step, use this universal geometric result and {\em presto!} the form of the optimal solution appears. How can we have a geometric result when probability distributions/densities and spectra naturally belong to the intersection of suitable cones or simplices with {\em $L^1$ spaces}? The reader might  look askance at this approach as, in general, in an infinite dimensional setting, $L^1$ spaces are not contained in $L^2$ spaces (one exception: absolutely summable sequences are also square summable). Hence, we simply don't have the Euclidean or Hilbert space geometry where orthogonality makes sense\footnote{This might well be the very  reason that our simple observation has not been made before in a countless number of papers on maximum entropy problems.}. However, in many important maximum entropy problems, the solution together with an appropriate function of it (inverse, logarithm, etc.) also belongs to a suitable $L^2$ space (when this is not the case, see Section \ref{gen-ban}, a more general Banach space result may be applied). Thus, as we show, there is nothing to loose formulating the problem over an appropriate Hilbert space possibly intersected with a cone or a simplex.

One might wonder at this point: What has this to do with the well known orthogonality principle of linear quadratic optimization? Right on! Theorem \ref{hilbertspace}, when applied to problems with quadratic criterion, yields well known results such as the orthogonality of the estimation error to the subspace generated by the available random variables in linear least-squares estimation. Thus, this orthogonality principle, a true {\em deus ex machina}, applies equally well to least-squares and entropic variational problems with linear constraints. Can this geometric result then be applied  to {\em any} optimization problem in  Hilbert space with linear constraints? Answer: No. The smoothness of the index functional is indispensable. For instance, the large and important class of  {\em compressed sensing} problems \cite{CRT,donoho,donoho2,CP,candesromberg,CR,romberg} features as criteria $l^1$-type norms which do not even admit directional derivatives (they only admit {\em one-sided} directional derivatives as they are convex).

The reader might be  doubtful by now: Don't the authors of this paper know about {\em information geometry, I-projections} and the like \cite{chentsov,csiszar0,topsoe,csiszar1,A,csiszar2,AN,KV,bathia,georgiou_scalar,YM,GEORGIOU_JIANG_NING,NJG}? We do and are savvy enough to know that this body of work is of central importance in Mathematical Statistics, Information Theory, Signal Processing, Identification and Control.  Our approach, however,  is different. Rather than viewing the solution itself of maximum entropy problems  as a projection in a suitable geometry and then developing a ``Pythagorean Theorem for I-divergences", our result involves usual orthogonality in Hilbert space (and the usual Pythagorean Theorem). Only that the orthogonality is a property of  the differential of the entropy functional which does not in general relate to an ``error". In particular, our geometry does not depend on the particular entropic criterion employed but only on the Hilbert space in which the primal variables live.

The paper is outlined as follows. In Section \ref{geometricresult}, we present our basic variational result. This is then applied in Sections \ref{matrix variational} and \ref{matricial functions} to various classical and more recent maximum Burg's entropy problems and in Section \ref{prior} to entropy problems with prior. In Section \ref{shannon}, we discuss  maximum entropy problems on a finite measure space. In Section \ref{reciprocal}, we develop a new application to block-circulant covariance matrix completion when an a priori estimate is available. Finally, in Section \ref{gen-ban}, we give a  generalization of our main result to Banach spaces.

\section{Maxima on surfaces}\label{geometricresult}

Let $G:\R^3\rightarrow\R$ be a continuously differentiable map and consider the surface (level set) ${\cal S}\subset\R^3$ determined by the equation
$$G(x)=c, \quad c\in\R.
$$
Since the derivative of $G$ in the direction of a vector $v$ tangent to the surface ${\cal S}$ must be zero, we get $\nabla G \cdot v=0$. It namely follows the well-known fact that the gradient $\nabla G(x_0), x_0\in{\cal S},$  is perpendicular to the plane tangent to the surface ${\cal S}$ at $x_0$. Let $F:\R^3\rightarrow\R$ be another smooth functional and suppose that we are interested in maximizing  $F$ over ${\cal S}$. By the chain rule, at a local maximum point $x_0$, $\nabla F$ must be orthogonal to every differentiable curve on ${\cal S}$ passing through $x_0$. We conclude that, at a maximum point $x_0$ the gradient of $F$ must also be perpendicular to the plane tangent to the surface ${\cal S}$ at $x_0$ and therefore aligned with $\nabla G (x_0)$, cf. e.g. \cite[pp. 101-109]{EDWARDS}. For instance, suppose we want to minimize $F(x,y,z)=x+y+z$ on the surface of the unit sphere $G(x,y,z)=x^2+y^2+z^2=1$. Since at a maximum point $\nabla F=(1,1,1)^T$ must be proportional to $\nabla G=2(x,y,z)^T$, we conclude that maxima have equal components. It follows that the unique maximum point is $X_M=(3^{-1/2},3^{-1/2},3^{-1/2})$ ($X_m=-X_M$ is the unique minimum point).

The purpose of this paper is to show that a suitable generalization of this basic result is sufficient to derive the form of the optimal solution in a variety of {\em maximum entropy problems}.

In maximum entropy problems, the map $G$ is  actually {\em linear} on a suitable vector space. Hence,  ${\cal S}=\ker G +c$ is an {\em affine} space, namely the translation of the subspace ${\cal V}=\ker G$. In this case, the geometric principle simply says that, at a maximum point, $\nabla F$  must be perpendicular to the subspace $\ker G$.
We apply this geometric result to a large class of {\em Burg-entropy} and {\em Shannon-entropy}  \cite{csiszar2,CsiszarMatus,CsiszarMatus2} variational problems encompassing, as special cases,  {\em Burg spectral estimation method} \cite{BURG1,BURG2}, {\em Dempster's covariance completion} \cite{Dempster-72} and {\em Gibbs-like variational principles} \cite{ellis}. In Burg's maximum entropy  problems, one maximizes the (integral of the) logarithm of a positive quantity, be it a  probability density, a spectrum or the determinant of a positive definite matrix, under linear constraints. The latter determine the affine space $\cal W$. Theorem \ref{hilbertspace} simply says that the  Fr\'echet differential of the entropy functional at a critical point must belong to the orthogonal complement of the subspace ${\cal V}$ of which the affine space $\cal W$ is a translation (coset). In the Burg's entropy case, this entails that the adjoint of the inverse of the solution must belong to ${\cal V}^\perp$. In the case of Shannon maximum entropy problems, the orthogonality condition concerns the logarithm of the solution. 

Classical results  can then be readily re-derived and generalized. For example, our result contains the key for the (considerable) recent generalizations developed in \cite{GEORGIOU_MAXIMUMENTROPY,Enqvist_Karlsson_covariance-Interpolation,FP,CFPP,FERRANTE_PAVON_ZORZI_ENHANCEMENT,FMP}. The case when a {\em prior estimate} is available is also covered by this geometric principle. In the Burg's case the entropic functional turns into a multivariate {\em Itakura-Saito divergence}  \cite{BASSEVILLE_DISTANCEMEASURES, DISTORTION}. In the Shannon case, entropy is replaced by the {\em Kullback-Leibler divergence} (relative entropy) \cite{kullback}). As an application, we show how our result can be used to extend the results of \cite{CFPP} to the case when a prior estimate of the circulant covariance is available. The latter problem deals with identifying of the parameters of a stationary  reciprocal process given the first covariance lags and an a priori covariance estimate.

Let ${\mathcal H}$ be a Hilbert space and let  $F:{\mathcal H}\rightarrow \R$ be a functional. We say that $F$ is {\em G\^{a}teaux-differentiable} at $h_0$ in direction  $v$ if the limit
$$F'(h_0;v):=\lim_{\epsilon\rightarrow 0}\frac{F(h_0+\epsilon v)-F(h_0)}{\epsilon}
$$
exists. In this case, $F'(h_0;v)$ is called the {\em directional derivative} of $F$ at $h_0$ in direction $v$. We say that $F$ is {\em Fr\'echet-differentiable} at $h_0$  if there exists an element $DF(h_0)$ in ${\mathcal H}$ such that
$$\lim_{\|h\|_{_{\mathcal H}}\rightarrow 0
}\frac{|F(h_0+h)-F(h_0)-\langle DF(h_0),h\rangle_{\mathcal H}|}{\|h\|_{\mathcal H}}=0.
$$
The element $DF(h_0)$ is called the {\em Fr\'echet differential} of $F$ at $h_0$. Fr\'echet differentiability is stronger than G\^{a}teaux differentiability. In fact, we have the following result \cite[p.50]{K}.
\begin{proposition} \label{gato-fresce}Let $F$ be Fr\'echet differentiable at $h_0$. Then, $DF(h_0)$ is unique and, for any $v\in{\mathcal H}$, $F$ is G\^{a}teaux differentiable at $h_0$ in direction $v$ and it holds
\begin{equation}\label{FGrel}
F'(h_0;v)=\langle DF(h_0),v\rangle_{\mathcal H}.
\end{equation}
Conversely,  when $F$ is G\^{a}teaux differentiable on an open set ${\mathcal U} \subseteq {\mathcal H}$ and its 
G\^{a}teaux derivative is linear and continuous at each point of ${\mathcal U}$ 
then $F$  is Fr\'echet differentiable in  ${\mathcal U}$.
Finally, when $F$ is convex, if it is G\^{a}teaux differentiable in all directions $v$ then it is Fr\'echet differentiable.
\end{proposition}

In some applications, we cannot expect that the functional be Fr\'echet differentiable at the point of interest. We may, however, have that a formula like (\ref{FGrel}) holds when $v$ varies over a subspace. More precisely, let ${\mathcal V}\subseteq {\mathcal H}$ be a (not necessarily closed) subspace and  $h\in{\mathcal H}$. Consider the corresponding {\em coset} ${\mathcal W}:= h+{\mathcal V}$ which is an {\em affine space} over ${\mathcal V}$. Observe that, for $w\in{\mathcal W}$ and $v\in{\mathcal V}$, $(w+\epsilon v)\in{\mathcal W}$ for all for all real $\epsilon$,, namely $w$ is an {\em internal point} of ${\mathcal W}$ in direction $v$. 
\begin{definition}We say that  $w_c$ is a {\em critical point} of $F$ over ${\mathcal W}=h+{\mathcal V}$ if $F'(w_c;v)=0$ for all $v\in{\mathcal V}$.
\end{definition}
\begin{theorem}\label{hilbertspace}Let ${\mathcal W}:= h+{\mathcal V}$ be an affine space. Assume that the functional $F$ is G\^{a}teaux-differentiable  at $w_c\in{\mathcal W}$ in any direction $v\in{\mathcal V}$ and that the G\^{a}teaux differential is given by the linear, continuous map $F'(w_c;v)=\langle D_{{\mathcal V}}F(w_c),v\rangle_{\mathcal H}$ where $D_{{\mathcal V}}F(w_c)\in {\mathcal H}$ .
Then $w_c$ is a critical point of $F$ over ${\mathcal W}$ if and only if $D_{{\mathcal V}}F(w_c)\in{\mathcal V}^\perp$. When  F is actually Fr\'echet differentiable at $w_c\in{\mathcal W}$, $w_c$ is  critical  if and only if $DF(w_c)\in{\mathcal V}^\perp$.\end{theorem}
\proof
$F'(w_c;v)=0$ for all $v\in{\mathcal V}$ if and only if $\langle D_{{\mathcal V}}F(w_c),v\rangle_{\mathcal H}=0, \forall v\in{\mathcal V}$.
\qed


\section{Matricial variational problems}\label{matrix variational}
\subsection{Geometric result}\label{matricialgeometricresult}

Let  ${\cal H}=\C^{n\times n}$ (or ${\cal H}=\R^{n\times n}$) be the space of $n\times n$ matrices  endowed with the inner product $\langle M_1,M_2\rangle:=\tr[M_1^* M_2]$, where $*$ denotes transposition plus conjugation (we write $M^{-\ast}$ for $(M^{-1})^\ast$).
The following result was established in \cite{FP}.

\begin{lemma}\label{chainr}
Let 
\eq\label{index-matr-abs}
F(M):=\log|\det[M]|.
\eeq
If $M$ is nonsingular then, for all $\delta M\in{\mathcal H}$
\eq
F'(M;\delta M)=\tr[M^{-1}\delta M]=\langle M^{-\ast}, \delta M\rangle.
\eeq
\end{lemma}
It now follows from Proposition \ref{gato-fresce} that $F$ is Fr\'echet differentiable in the open set of non-singular matrices and 
\eq
\label{deffre}
DF(M)=M^{-\ast}.
\eeq

We are interested in extremizing (\ref{index-matr-abs}) over an {\em affine space}, namely a coset of the form ${\mathcal W}=A+{\cal V}$, where $A\in {\mathcal H}$ and ${\mathcal V}$ is a subspace of ${\mathcal H}$.
\begin{theorem}\label{theorem1}
Let ${\cal W}=A+{\cal V}$ be an affine space. Then a nonsingular matrix $M_c\in{\cal W}$ extremizes $F(M)=\log|\det[M]|$ over ${\cal W}$ if and only if $M_c^{-*}\in{\cal V}^\perp$.
\end{theorem}
\proof
 Let $M_c\in {\cal W}$ be non-singular. By (\ref{deffre}), we have $DF(M_c)=M_c^{-\ast}$.
The conclusion now follows from Theorem \ref{hilbertspace}.
 \qed
 \subsection{Dempster's covariance selection}\label{dempstermaxent}
In various applications, index (\ref{index-matr-abs}) must be extremized (or rather maximized) on the intersection between an affine space ${\cal W}$ and a convex cone. A typical example is that of the cone of positive semidefinite matrices. This is the case considered by Dempster in the seminal paper \cite{Dempster-72} where  a general strategy for completing a partially specified covariance matrix was introduced. 
We now show that Theorem \ref{theorem1} provides a geometrical interpretation of one of the key features of Dempster's result. 
To see this, consider the Dempster's problem with the same notation as in Subsection \ref{dempster}. Let  ${\cal W}$ be the affine space of symmetric matrices having  elements $\{\sigma_{ij};1\le i\le j\le n, (i,j)\in \bar{{\cal I}}\}$. Notice that ${\cal W}$ is affine over the subspace ${\cal V}$ of symmetric matrices having zeros in the positions $\bar{{\cal I}}$. Observe next that the solution $\Sigma$ is  constrained to be in the intersection between ${\cal W}$ and the convex cone of positive definite matrices.  On this set, maximizing the index (\ref{index-matr-abs}) or the  entropy (\ref{entropy}) is equivalent.  Thus, the two criteria yield the same solution. Moreover, $(i,i)\in\bar{{\cal I}}$ for all $i=1\dots n$, i.e. $[\Sigma]_{ii}$ are all fixed so that
$|[\Sigma]_{ij}|\leq \sqrt{[\Sigma]_{ii}[\Sigma]_{jj}}$ and hence the feasible set is bounded. Finally, as $\Sigma$ tends to be singular, i.e. it approaches the boundary of the cone, $H(p)$ tends to $-\infty$ which implies that the solution can be searched among positive definite matrices.
Thus, under the feasibility assumption, the optimal solution exists and lies in the {\em interior} of the cone. We can then repeat locally the argument of Theorem \ref{theorem1} to conclude that the maximum entropy completion $\Sigma_c$ is such that $\Sigma_c^{-1}\in{\cal V}^\perp$. Finally, observe that ${\cal V}^\perp$ is the space of matrices having zeros in ${\cal I}$, the complement  of $\bar{{\cal I}}$. Indeed, let  $e_i$ denote the $i$-th canonical vector in $\R^n$ and observe that for $(i,j)\in {\cal I}$, the rank one matrix $e_ie_j\tp$ belongs to ${\cal V}$. If $M\in{\cal V}^\perp$, we must have
$$0=\tr[(e_ie_j\tp)\tp M]=
\tr[e_je_i\tp M]=e_i\tp Me_j=[M]_{ij},\quad\forall (i,j)\in{\cal I}.$$
Thus, the maximum entropy completion $\Sigma_c$ is a Dempster's completion.
\subsection{General matrix completion}
In \cite{FP}, Dempster's completions where shown to solve suitable entropy-like variational problems for general nonsingular matrices\footnote{Actually, the case of  full-rank rectangular matrices, with the Moore-Penrose pseudoinverse in place of the inverse, was also treated in \cite{FP}.}. Again, the form of the extremal completions (no uniqueness is there guaranteed) when they exist is provided by Theorem \ref{theorem1}.

%
 
 \section{Matricial functions}\label{matricial functions}
 \subsection{The orthogonality result}\label{orth}
 
 Consider now the Hilbert space  ${\mathcal H}$ of square integrable functions defined on the unit circle $\T$ and taking values in the space of $m \times m$ Hermitian matrices.
We denote by $\Hermitian_n$ the $n^2$-dimensional, {\em real} vector space of
Hermitian matrices of dimension $n \times n$. Hence, ${\mathcal H}=L^2(\T,\Hermitian_m)$ with scalar product 
 $$\langle \Phi,\Psi\rangle_{\cal H}:=\frac{1}{2\pi}\int^\pi_{-\pi}\tr\left[\Phi(e^{j\vartheta})\Psi(e^{j\vartheta})\right]d\vartheta.
 $$
Consider  the functional
\eq\label{index-func-abs}
F(\Phi)=\frac{1}{2\pi}\int^\pi_{-\pi}\log|\det[\Phi(e^{j\vartheta})]|d\vartheta.
\eeq
\begin{lemma}\label{chainr2} Suppose $\Phi\in L^\infty(\T,\Hermitian_n)$ is coercive\footnote{$\Phi$ is called coercive if $\exists\,\alpha>0$ s.t.  $\Phi(e^{{\imunit}\vartheta})-\alpha I_m$ is a.e. positive definite on $\T$.}. 
Then,  for any $\delta \Phi\in L^\infty(\T,\Hermitian_n)$ the directional derivative of (\ref{index-func-abs}) exists and is given by the linear map
\eq
F'(\Phi;\delta \Phi)=\frac{1}{2\pi}\int^\pi_{-\pi}\tr[\Phi^{-1}(e^{j\vartheta})\delta \Phi(e^{j\vartheta})]d\vartheta=\langle \Phi^{-1}, \delta\Phi\rangle_{\mathcal H}.
\eeq 
\end{lemma}
\proof 
Observe that, for $\delta \Phi\in L^\infty(\T,\Hermitian_n)$ and $|\varepsilon|$ sufficiently small, $\Phi(e^{j\vartheta})+\varepsilon\delta \Phi(e^{j\vartheta})$ is a.e. positive definite.  After bringing the derivative under the integral sign, we can  use Lemma \ref{chainr} for almost all $\vartheta$.
\qed
Let ${\cal W}=A+{\cal V}$  be an affine space in $L^\infty(\T,\Hermitian_n)$, namely $A\in L^\infty(\T,\Hermitian_n)$ and ${\mathcal V}$ is a subspace of $L^\infty(\T,\Hermitian_n)$. 

Then Theorem \ref{hilbertspace} yields:
 \begin{theorem}\label{theorem2}
Let ${\cal W}=A+{\cal V}$ be as above  and $\Phi_c(e^{j\vartheta})\in {\cal W}$ be coercive. Then,
if $\Phi_c$ is a critical point of
(\ref{index-func-abs}) over ${\cal W}$, we have   $\Phi_c^{-1}\in{\mathcal V}^\perp$\footnote{For a not necessarily closed subspace ${\mathcal V}$ of $L^2(\T,\Hermitian_n)$, the {\em orthogonal complement} ${\mathcal V}^\perp$ is the closed subspace of  $u\in L^2(\T,\Hermitian_n)$ such that $\langle u,v\rangle_{L^2}=0, \forall v\in{\mathcal V}$.}
.
\end{theorem}
\proof By Lemma \ref{chainr2}, if such a $\Phi_c$ extremizes (\ref{index-func-abs}), then $\langle \Phi^{-1}_c, v\rangle_{\mathcal H}=0$ for all $v\in {\cal V}$. Namely, $\Phi_c^{-1}\in{\mathcal V}^\perp$.
\qed
 \subsection{Burg's maximum entropy covariance extension}\label{burgME}

In his seminal work \cite{BURG1,BURG2}, Burg introduced a spectral estimation method based on the maximization of entropy which is widely used in signal processing.  We now show that Theorem \ref{theorem2}  provides a most transparent reason why the solution has to be an AR process.
 Consider a discrete-time { Gaussian} process $\{y_k;\,k\in\Z\}$ taking values in $\R^m$. Let  $Y_{[-n,n]}$ be the random vector obtained by considering the window $y_{-n},y_{-n+1},\cdots,y_0,\cdots,$ $y_{n-1},y_n$, and let $p_{Y_{[-n,n]}}$ denote the corresponding joint density. The {\em (differential) entropy rate} of $y$ is defined by
\begin{equation}\label{ER}
{\h{y} :=}\lim_{n\rightarrow\infty}\frac{1}{2n+1}H(p_{Y_{[-n,n]}}),
\end{equation}
if the limit exists, where $H(p_{Y_{[-n,n]}})$ denotes the  entropy of the density of the random vector $Y_{[-n,n]}$, cf. (\ref{entropy}).
In \cite{KOL56}, Kolmogorov established the following important result.
\begin{theorem} \label{Kolmogorov}Let $y=\{y_k;\,k\in\Z\}$ be a $\R^m$-valued, zero-mean, Gaussian, stationary, purely nondeterministic of full rank  process with spectral density $\Phi_y$.
Then
\begin{equation}\label{entropyconnection}
\h{y}=\frac{m}{2}\log (2\pi e)+\frac{1}{4\pi}\int_{-\pi}^\pi\log\det\Phi_y(e^{{\imunit}\vartheta})d\vartheta.
\end{equation}
\end{theorem}
As is well-known, there is also a fundamental connection between the quantity appearing in
(\ref{entropyconnection}) and the optimal one-step-ahead predictor: The multivariate Szeg\"{o}-Kolmogorov formula reads
\begin{equation}\label{SK}
\det R=\exp\left\{\frac{1}{2\pi}\int_{-\pi}^\pi\log\det\Phi_y(e^{{\imunit}\vartheta})d\vartheta\right\},
\end{equation}
where $R$ is the error covariance matrix corresponding to the optimal predictor. Consider now the multivariate covariance extension problem. Let  $C_k, k=0,1,\ldots,n-1$ of dimension $m\times m$ be some estimated  {\em covariance lags} of  an unknown stationary process $y$. Then {\em Burg's problem} consists in finding a stationary process $y$ with spectral density $\Phi_y$ which maximizes the index
\eq
\label{func-entro}
F(\Phi_y)=\frac{1}{2\pi}\int_{-\pi}^\pi\log\det\Phi_y(e^{{\imunit}\vartheta})d\vartheta.
\eeq
among all spectral densities having as first $n$ Fourier coefficients  $C_k, k=0,1,\ldots,n-1$. In view of Kolmogorov's result (\ref{entropyconnection}), maximizing the entropy rate of a stationary Gaussian process is equivalent to maximizing the integral of $\log\det\Phi_y$.\footnote{Actually, the solution to this problem maximizes the entropy rate in the larger class of second-order processes  \cite {CHOI_COVER}.} Assume that the block-Toeplitz matrix $\Sigma_n$
\eq
\label{TOEPLITZ}
\Sigma_n = \bmat{cccc}
     C_0       &C_1&\cdots&C_{n-1}\\
     C_1^\ast &C_0&\cdots&C_{n-2}\\
     \vdots&\vdots& \ddots&\vdots\\
     C_{n-1}^\ast&C_{n-2}^\ast&\dots&C_0
     \emat
\eeq
is positive definite. Then \cite{GRENANDER_SZEGO} there are infinitely many spectra having the prescribed Fourier coefficients. 

Consider now the matrix pseudo-polynomial $P(e^{{\imunit}\vartheta})=\sum_{k=-n+1}^{n-1}C_ke^{-{\imunit}\vartheta k}$, with $C_{-k}:=C^\ast_k$, and define the subspace ${\cal V}_n$ of $L^\infty(\T,\Hermitian_n)$ of functions whose Fourier coefficients
$R_i$ vanish for all $i=-n+1,\dots n-1$ and obey to the symmetry constraint
$R_i=R^\ast_{-i}$. Then the constraint in Burg's problem can be expressed as $\Phi\in{\cal W} \cap {\cal S}$, where the affine space ${\cal W}$ is defined by
$${\cal W}=P+{\cal V}_n$$
and ${\cal S}$ is the convex cone of bounded, coercive spectral densities.
On   ${\cal S}$, (\ref{index-func-abs}) and (\ref{func-entro}) coincide,   and
$F$ is strictly concave. Thus, an extremizer $\Phi_c$ is actually a maximum point.
By Theorem \ref{theorem2}, this maximum point  $\Phi_c$ is such that $\Phi_c^{-1}\in\bar{{\cal V}}_n^\perp$. Observe now that $\bar{{\cal V}}_n^\perp$ is given by the matricial polynomials of the form
$$ Q(e^{{\imunit}\vartheta})=\sum_{k=-n+1}^{n-1}A_ke^{-{\imunit}\vartheta k},\quad A_{-k}=A_k^\ast. 
$$
We conclude that the optimal spectrum has the form
\begin{equation}\label{equzburg}\Phi_c(e^{{\imunit}\vartheta})=\left[\sum_{k=-n+1}^{n-1}A^\circ_ke^{-{\imunit}\vartheta k}\right]^{-1},\quad A_{-k}^\circ=(A_k^\circ)^\ast
\end{equation}
for some matrices $A^\circ_k, k=-n+1,\ldots,0,\ldots,n-1$ which permit to satisfy the constraints on the first $n$ coefficients. Thus, the solution process is an AR process. 
If only some  of the $C_k, k=0,1,\ldots,n-1$ are available, the classical approach to the problem requires a certain effort and some {\em ad hoc} reasoning  to get the solution form.
Theorem \ref{theorem2}, on the contrary, yields immediately   that in (\ref{equzburg}) $A_k^\circ=0$ for all $k$ corresponding to missing $C_k$'s.
    
\subsection{A more general moment problem}\label{moregeneral}
We consider next a generalization of Burg's problem studied by Byrnes, Georgiou and Lindquist  and co-workers \cite{BGuL,BGL2,BYRNES_GUSEV_LINDQUIST_FROMFINITE,G4,G5',GEORGIOU_LINDQUIST_KULLBACKLEIBLER,G5,L,GEORGIOU_LINDQUIST_CONVEXARMA} in the frame of generalized moment problems. In their broad research effort, having applications, besides spectral estimation, to robust control problems,  elements of a parametric family of rational spectral  densities were recognized from the start \cite{BGuL,BGL1} to be critical points of logarithmic entropy-like functionals.

Consider a transfer function 
\begin{equation} \label{filterbank} G(z)=(zI-A)^{-1}B, \;\; A\in\C^{n
\times n} , B\in\C^{n \times m},\;\; n> m,
\end{equation}
 where $A$ has
all its eigenvalues in the open unit disk, $B$ has full column
rank, and $(A,B)$ is a reachable pair\footnote{A pair $(A,B)$ is called reachable in Systems Theory \cite{KFA} if the matrix $[B\mid AB\mid\dots\mid A^{n-1}B]$ has full row rank.}. Suppose $G(z)$ models a
bank of filters fed by a wide sense stationary, purely
nondeterministic, $\C^m$-valued process $y$:
\vspace{1.2cm}
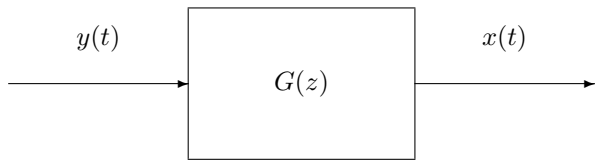
\begin{center}
\setlength{\unitlength}{.2cm}
\begin{picture}(43,8)(0,5)
\put(6,13){\makebox(0,0){$y(t)$}}
\put(33,13){\makebox(0,0){$x(t)$}}
\put(0,10){\vector(1,0){12}}
\put(27,10){\vector(1,0){12}}
\put(12,5){\framebox(15,10){$G(z)$}}
\end{picture}
\end{center} 
\vspace{0.5cm}
\noindent
Let $x$ be the
$n$-dimensional stationary output process
\begin{equation}
\label{statedyn}x_{k+1}=Ax_k+By_k, \;\; k\in \Z.
\end{equation}
We denote by $\Sigma$ the covariance of  $x_k$. The spectrum $\Phi$ must then satisfy the following {\em moment constraint}
\begin{equation}\label{moment}
\frac{1}{2\pi}\int_{-\pi}^{\pi} G(e^{{\imunit}\vartheta})\Phi(e^{{\imunit}\vartheta})G^\ast(e^{{\imunit}\vartheta})d\vartheta=\Sigma.
\end{equation}
 As in \cite{BYRNES_GUSEV_LINDQUIST_FROMFINITE,GEORGIOU_LINDQUIST_KULLBACKLEIBLER,G5,FERRANTE_PAVON_RAMPONI_HELLINGERVSKULLBACK,FMP}, we now consider the problem of determining  spectral densities
$\Phi$ satisfying (\ref{moment}) for a given $\Sigma>0$.
The covariance extension  is a special case of this problem corresponding to
$G(z):=[z^{-n}I\mid z^{-n+1}I\mid\dots\mid z^{-1}I]\tp$ and $\Sigma$ equal to the Toeplitz matrix in (\ref{TOEPLITZ}). More details on this fact may be found in \cite{GEORGIOU_LINDQUIST_KULLBACKLEIBLER} where other classical problems are shown to be special cases of the above.
The most important of these problems is the celebrated {\em Nevanlinna-Pick interpolation problem} of fundamental importance in various $H^\infty$ control problems \cite{Doyle,BYRNES_GEORGIOU_LINDQUIST_GENERALIZEDCRITERION,BLOMQVIST_LINDQUIST_NAGAMUNE_MATRIXVALUED,GEORGIOU_LINDQUIST_REMARKSDESIGN}. In the scalar case, the simplest version of the Nevanlinna-Pick problem\footnote{See e.g. \cite{BLOMQVIST_LINDQUIST_NAGAMUNE_MATRIXVALUED} for the general multivariable case.} consists in finding a {\em positive-real function} $Z(z)$\footnote{The concept of positive-real function was introduced by Cauer and Brune in 1930 in Network Theory as passive networks, such as RLC circuits, have impedance functions that are positive real.}, namely a function analytic in $\{|z|<1\}$ and having nonnegative real part there, that satisfies the following interpolation conditions:
$$Z(p_i)=w_i,\quad  i=1,\ldots,n,$$
where $p_i$ are given distinct points in the open unit disc and $w_i$ are given complex values. 
This problem becomes a special case of $(\ref{moment})$ with the following prescriptions for $Z$, $G$ and $\Sigma$. The positive-real function $Z$ is related to the sought spectral density $\Phi$ by
$$\Phi(z)=Z(z)+\overline{Z(\bar{z}^{-1})},$$ 
the $k$-th component of $G(z)$ is
$$G_k(z)=\frac{1}{z-p_k},$$
and the matrix $\Sigma$ is the {\em Pick
matrix} with
elements
$$\Sigma_{i,j}=\frac{w_i+\bar{w}_j}{1-p_i\bar{p}_j}.$$
A possible choice for $A$ and $B$ in (\ref{filterbank}) so that $G_k(z)=\frac{1}{z-p_k}$ is then
$$
A=\bmat{ccccc}
p_1&0&0&\dots&0\\
0&p_2&0&\dots&0\\
\vdots&\vdots&
&\ddots&\vdots\\
           0&0&0&\dots&0\\
           0&0&0&\dots&p_n
           \emat,\quad
B=\bmat{c}1\\1\\\vdots\\1\\1\emat.
$$
           
\vspace {0.3cm}
\noindent
Notice that the complex numbers $w_i$ may be recovered from a solution  $\Phi$ through
$$w_k=\frac{1}{4\pi}\int_{-\pi}^{\pi}\frac{e^{-i\omega}+p_k}{e^{-i\omega}-p_k}
\Phi(e^{i\omega})d\omega,\quad
k=1,2,\dots, n.
$$

We now show how to treat this problem in our geometric framework. Let, as before, ${\mathcal H}=L^2(\T,\Hermitian_m)$.
Consider now the linear operator 
\eqn \Gamma:  L^\infty(\T,\Hermitian_m) & \rightarrow & \Hermitian_n \nn,\\
\Phi & \mapsto & \frac{1}{2\pi}\int_{-\pi}^{\pi} G(e^{{\imunit}\vartheta})\Phi(e^{{\imunit}\vartheta})G^\ast(e^{{\imunit}\vartheta})d\vartheta.\eeqn  It follows that for the constraint (\ref{moment}) to be feasible, $\Sigma$ must belong
to the linear space 
\eq \label{rangegamma} \Rgamma:=\left\{ M \in \Hermitian_n| \exists
\Phi \in L^\infty(\T,\Hermitian_m), \frac{1}{2\pi}\int_{-\pi}^{\pi} G\Phi G^\ast d\vartheta=M \right\}.\eeq
Consider now the following generalization of Burg's problem: Maximize the entropy index (\ref{func-entro}) subject to (\ref{moment}) where $\Sigma$ is assumed to be positive definite. Suppose that (\ref{moment}) is feasible, namely there exists a spectral density $\Phi_0\in L^\infty(\T,\Hermitian_m)$ satisfying this constraint. Then, the family ${\cal W}$ of hermitian-valued functions satisfying (\ref{moment}) may be expressed as
$${\cal W}=\Phi_0+{\cal V},$$
where ${\cal V}=\{\Phi \in L^\infty(\T,\Hermitian_m)|\int G\Phi G^\ast=0\}$. In other words, ${\cal V}=\ker \Gamma$.   The constraint in the generalized Burg problem can be expressed as $\Phi\in{\cal W}\cap {\cal S}$, where ${\cal S}$ is the convex cone of bounded, coercive spectral densities.
Since 
\eqn
\nn
\langle \int_{-\pi}^{\pi} G\Phi G^\ast \frac{d\vartheta}{2\pi},M\rangle_{\Hermitian_n}&:=&\tr\left[\int_{-\pi}^{\pi} G \Phi G^\ast \frac{d\vartheta}{2\pi} M\right]=\tr\left[\int_{-\pi}^{\pi}  \Phi G^\ast MG \frac{d\vartheta}{2\pi} \right]\\
\nn
&=&\langle \Phi,G^\ast MG\rangle_{\mathcal H},
\eeqn
we have that the {\em adjoint} of $\Gamma$, mapping $\Hermitian_n$ to $L^\infty(\T,\Hermitian_m)$, is given by
\begin{equation} \Gamma^\ast:  M  \mapsto  G^\ast MG.
\end{equation} 
In particular, $\Range \Gamma^\ast=\{\Phi=G^\ast MG, M\in \Hermitian_n\}\subset  \Cthm\nn$ the continuous Hermitian-valued functions on the unit circle.
Since $\Range\Gamma^\ast$ is finite-dimensional, it is necessarily closed and we have
\begin{equation}\label{ortcompl}{\cal V}^\perp=\left[\ker \Gamma\right]^\perp=\Range \Gamma^\ast=\{\Phi=G^\ast MG, M\in \Hermitian_n\}.
\end{equation}
By Theorem \ref{theorem2}, the maximum point  $\Phi_c$ is such that $\Phi_c^{-1}\in{\cal V}^\perp$. Hence, the optimal spectrum has the form
\begin{equation}\label{optspectrum}
\Phi_c(e^{j\vartheta})=\left[G(e^{j\vartheta})^\ast \Lambda_c G(e^{j\vartheta})\right]^{-1},
\end{equation}
for some Hermitian $\Lambda_c$ such that $G(e^{j\vartheta})^\ast \Lambda_c G(e^{j\vartheta})>0$ on $\T$ and the constraint (\ref{moment}) is satisfied, namely
$$\int_{-\pi}^{\pi} G\left[G^\ast \Lambda_c G\right]^{-1}G^\ast\frac{d\vartheta}{2\pi}=\Sigma.
$$
Indeed, Georgiou showed in \cite{GEORGIOU_MAXIMUMENTROPY} that the unique solution of the generalized Burg problem has the form (\ref{optspectrum}) with
\begin{equation}
\Lambda_c=\Sigma^{-1}B\left(B^\ast\Sigma^{-1}B\right)^{-1}B^\ast\Sigma^{-1}.
\end{equation}

 \section{Variational entropy problems with ``prior"}\label{prior}
 \subsection{Matricial problems}
 Consider now the same set up as in Section \ref{matrix variational}, where a ``prior" nonsingular estimate $N$ of the matrix $M$ is available. Rather than extremizing (maximizing) (\ref{index-matr-abs}), we now consider the problem of finding a matrix belonging to the given affine set ${\cal W}$ and which 
extremizes the index
\begin{equation}
\label{index-prior}
F_N(M):=\log|\det[N]|-\log|\det[M]|+\tr \left(N^{-1}M\right)
\end{equation} 
(see below for insights and motivation for this choice).
Lemma \ref{chainr} now becomes:
\begin{lemma}\label{chainr3}
Let 
$F_N(M)$ be given by (\ref{index-prior}).
If $M$ is nonsingular then for any $\delta M\in{\mathcal H}=\C^{n\times n}$,
\eq
F_N'(M;\delta M)=\tr[\left(-M^{-1}+N^{-1}\right)\delta M],
\eeq
and  $DF_N(M)=-M^{-1}+N^{-1}$. 
\end{lemma}
By Theorem \ref{hilbertspace}, we get:
\begin{theorem}\label{theoremmatrixprior}
Let ${\cal W}=A+{\cal V}$ be an affine set in ${\mathcal H}=\C^{n\times n}$. Let $N$ be a nonsingular matrix in ${\mathcal H}$. Then the nonsingular matrix $M_c\in{\cal W}$ extremizes (\ref{index-prior}) over ${\cal W}$ if and only if  $\left(M_c^{-*}-N^{-*}\right)\in{\cal V}^\perp$.
\end{theorem}
\begin{remark} Notice that the optimality condition may be expressed as
$$\P^{\cal V}M_c^{-*}=\P^{\cal V}N^{-*},
$$
where $\P^{\cal V}$ is the othogonal projection onto ${\cal V}$. Also notice that in the case when $N^{-*}\in{\cal V}^\perp$, the solution $M_c$ of the problem without prior of Subsection \ref{matricialgeometricresult} solves also this problem.
\end{remark}
\vspace{4mm}

 In order to motivate the choice (\ref{index-prior}), we first recall a few basic facts on entropy for Gaussian random random vectors and processes that may be found e.g. in \cite{PINSKER_INFORMATION,IHARA_INFORMATION,COVER_THOMAS}. The {\em relative entropy} or { {\em Kullback-Leibler}} pseudo-distance or {\em divergence} between two probability densities $p$ and $q$, with the support of $p$ contained in the support of $q$, is defined by
  \begin{equation}\label{kullback-leibler}
  \D(p\|q){ :=}\int_{\R^n} p(x)\log\frac{p(x)}{q(x)}dx,
\end{equation}
see e.g \cite{COVER_THOMAS}.
In the case of two zero-mean Gaussian { densities} $p$  and $q$ with { positive definite} covariance  matrices  $M$ and $N$, respectively, the relative entropy is given by:
\begin{equation}\label{divgauss}\D(p\|q)=\frac{1}{2}\left[ \log\det
(M^{-1}N)+\tr(N^{-1}M)-n \right].
\end{equation}
Hence, when $N$ and $M$ are positive definite, minimizing index (\ref{index-prior}) is indeed equivalent to minimizing the Kullback-Leibler divergence between two Gaussian random vectors which is one of the  central problems in statistical modeling. Indeed, as is well-known, (\ref{divgauss}), originates from {\em maximum likelihood} considerations, cf. e.g. \cite[Section II]{BLW}.
 An important application of this result is the {\em estimation of a structured covariance matrix}. In this class of problems, one seeks a covariance matrix $\Sigma$ that, besides being symmetric and positive definite, ejoys further properties such as being Toeplitz, circulant, etc.. The covariance estimated from the data $\hat{\Sigma}$ usually fails to have the prescribed structure. Hence, the problem arises to find $\hat{\Sigma}_c$ with the further properties which is as close as possible to $\hat{\Sigma}$, see
 \cite{BLW,Georgiou-LinquistFest,FERRANTE_PAVON_ZORZI_ENHANCEMENT,NJG} for more details and applications. 
This {\em static} problem has an important application as an ancillary problem also in the setting described in Section \ref{moregeneral}. Indeed, in the setting of Section \ref{moregeneral}, the state covariance $\Sigma$ of constraint (\ref{moment})
 is assumed to be given.
On the contrary, in practical situations, it must be estimated from the available data i.e. a finite sample of the unknown stochastic process $y$. More explicitly, the estimate of $\Sigma$ can be obtained as follows:
\begin{itemize}
    \item The filter $G(z)$ is fed by the $m$-dimensional data $\{y_i\}_{i=1}^N$
    and we collect the $n$-dimensional output data $\{x_i\}_{i=1}^N$.
    \item We compute the sample covariance estimate $\hat{\Sigma}$ of $\Sigma$ in the usual way
    \eq \label{covcampionaria}\hat{\Sigma}:=\frac{1}{N}\sum_{i=1}^N x_ix_i^*.\eeq
\end{itemize}
Notice that $\hat{\Sigma}\in \Hermitian_n$ and $\hat{\Sigma}\ge
0$. Nevertheless, in general, $\hat{\Sigma}$ does not belong to the range of the operator $\Gamma$ given by (\ref{rangegamma}) so that the problem of Subsection \ref{moregeneral}  is unfeasible. Before we try to solve the generalized Burg problem of  Subsection \ref{moregeneral} we then need to approximate $\hat{\Sigma}$ with a suitable covariance matrix $\hat{\Sigma}_c$ which belongs to $ \Rgamma$.
If we take $F_{\hat{\Sigma}}(\cdot)$ as in (\ref{index-prior}) as distance index,  we have to minimize
$F_{\hat{\Sigma}}$ over the set of symmetric, positive definite matrices belonging to the range of   $\Gamma$.
This problem has been considered and solved in  \cite{FERRANTE_PAVON_ZORZI_ENHANCEMENT}.
In particular, it was shown in \cite[Proposition 3.2]{FERRANTE_PAVON_ZORZI_ENHANCEMENT} that, given the matrices $A$ and $B$ as in Section \ref{moregeneral},
the range of   $\Gamma$ may be characterized as
 \eq
{\cal V}= \{\Sigma:\ (I-\pjb)(\Sigma-A\Sigma
A^*)(I-\pjb)=0 \},
\eeq
with $\pjb$ being the orthogonal projection onto $\im(B)$, so 
that it is easy to see that 
\eq
{\cal V}^\perp= \{\Delta=(I-\Pi_B)\Lambda(I-\Pi_B) - A^*(I-\Pi_B)\Lambda (I-\Pi_B) A:\ \Lambda\in \Hermitian_n\}.
\eeq
Then,  Theorem \ref{theoremmatrixprior} can be used  to get in a straightforward manner the form of the optimal $\hat{\Sigma}_c$  \cite[Section IV]{FERRANTE_PAVON_ZORZI_ENHANCEMENT}:
\eq
\hat{\Sigma}_c=\left(\hat{\Sigma}^{-1}+(I-\Pi_B)\Lambda(I-\Pi_B) - A^*(I-\Pi_B)\Lambda (I-\Pi_B) A\right)^{-1},\ \Lambda\in \Hermitian_n.
\eeq

\subsection{Matricial functions problems with ``prior"}

As much as Theorem \ref{theorem1}, also Theorem \ref{theoremmatrixprior} may be generalized to the case when 
${\mathcal H}=L^2(\T,\Hermitian_n)$.
In this setting, we consider $\Phi\in L^\infty(\T,\Hermitian_n)$ coercive and a given ``prior" 
$\Psi$ also essentially bounded and coercive.
The index to be extremized is
\eq\label{index-fun-pri}
F(\Phi,\Psi)=\frac{1}{2\pi}\int_{-\pi}^\pi\left\{\log(\det \Psi)-\log(\det\Phi)+\tr\left[\Psi^{-1}\Phi\right]\right\}d\vartheta.
\eeq
Motivation for considering this index will be provided after the statement of the next result. A straightforward generalization of Lemma \ref{chainr3} and Theorem \ref{hilbertspace} now give a result germane to  Theorem \ref{theoremmatrixprior}:
\begin{theorem}\label{theorem4}
Let ${\mathcal H}$ be as before $L^2(\T,\Hermitian_n)$ and let ${\cal W}=A+{\cal V}$ be an affine set in $L^\infty(\T,\Hermitian_n)$ and $\Phi_c(e^{j\vartheta})\in {\cal W}$ be coercive. Then
$\Phi_c$ extremizes 
(\ref{index-fun-pri}) over ${\cal W}$  if and only if   $\left(\Phi_c^{-1}-\Psi^{-1}\right)\in {\cal V}^\perp$.
\end{theorem}

\vspace{3mm}
\noindent
Again, the optimality condition may be written as
$$\P^{\cal V}\Phi_c^{-1}=\P^{\cal V}\Psi^{-1}.
$$
In the case when $\Psi^{-1}\in{\cal V}^\perp$, the solution $\Phi_c$ of the problem without prior of Subsection \ref{orth} solves also this problem. For instance, in the case of Burg's problem of Subsection \ref{burgME} supplemented with a prior being an AR process of order $\le n$, the solution is the same as without prior (\ref{equzburg}). It is namely the maximum entropy solution.

To provide some motivation and insight for index (\ref{index-fun-pri}), we
consider  two zero-mean, Gaussian, stationary, purely nondeterministic  processes $y=\{y_k;\,k\in\Z\}$ and $z=\{z_k;\,k\in\Z\}$ taking values in $\R^m$.
We consider the {\em relative entropy rate}  $\dr{y}{z}$ between $y$ and $z$ defined as
\begin{equation}\label{DR}
{\dr{y}{z} :=} \lim_{n\rightarrow\infty}\frac{1}{2n+1}\D(p_{Y_{[-n,n]}}\|p_{Z_{[-n,n]}})
\mbox{ if the limit exists}
\end{equation}
where $p_{Y_{[-n,n]}}$ and $p_{Z_{[-n,n]}}$ are the densities of the random vectors obtained from $y$ and $z$, respectively, by considering the ``windows" from time $-n$ to time $n$. 
Following in his mentor's footsteps, the great information theorist M. Pinsker \cite{PINSKER_INFORMATION} proved the following important result (see also \cite{SvS,IHARA_INFORMATION,LPBook}):
\begin{theorem} Let $y=\{y_k;\,k\in\mathbb{Z}\}$ and $z=\{z_k;\,k\in\mathbb{Z}\}$ be $\R^m$-valued, zero-mean, Gaussian, stationary, purely  nondeterministic processes with spectral density functions $\Phi_y$ and $\Phi_z$, respectively.
Assume, moreover, that  at least one of the following conditions is satisfied:
\begin{enumerate}
\item $\Phi_y\Phi_z^{-1}$ is bounded;
\item $\Phi_y\in\mathrm{L}^{2}\left(-\pi,\pi\right)$ and $\Phi_z$ is {\em coercive}.
\end{enumerate}
Then
\begin{equation}\label{relentropyconnection}
\dr{y}{z}=\frac{1}{4\pi}\int_{-\pi}^\pi\left\{\log\det \left(\Phi_y^{-1}(e^{{\imunit}\vartheta})\Phi_z(e^{{\imunit}\vartheta})\right)+\tr\left[\Phi_z^{-1}(e^{{\imunit}\vartheta})\left(\Phi_y(e^{{\imunit}\vartheta})-\Phi_z(e^{{\imunit}\vartheta})\right)\right]\right\}d\vartheta.
\end{equation}
\end{theorem}

The index (\ref{relentropyconnection}) has the form of a multivariate {\em Itakura-Saito divergence} of speech processing \cite{DISTORTION,BASSEVILLE_DISTANCEMEASURES} and is basically the same as (\ref{index-fun-pri}). Indeed, one of the main results of \cite{FMP} is based on the minimization of (\ref{index-fun-pri}) where 
$\Psi$ is a given ``prior" spectral density and $\Phi$ must belong to the intersection between the cone ${\cal S}$ of positive definite spectral densities and the affine set ${\cal W}$ of the solutions of the moment problem (\ref{moment}), for given $G$ and $\Sigma$ as in Subsection \ref{moregeneral}. Since the constraint is as before, so are the spaces ${\cal W}$ and ${\cal V}$. In particular, we have 
$${\cal V}^\perp=\{\Phi=G^\ast MG, M\in \Hermitian_n\}.
$$
By Theorem \ref{theorem4}, we get the form of the optimal spectrum derived in \cite{FMP}
$$\Phi_c=\left[\Psi^{-1}+G^\ast \Lambda_c G\right]^{-1},\quad \Lambda_c\in\Hermitian_n,
$$
where $\Lambda_c$ permits to satisfy (\ref{moment}).
\subsection{Kullback-Leibler approximation of spectral densities}
Consider the same set up as in Subsection \ref{moregeneral} in the scalar case ($m=1$) when an a priori estimate of the spectrum $\Psi$ is available.  The latter is assumed to be essentially bounded and coercive. In \cite{GEORGIOU_LINDQUIST_KULLBACKLEIBLER}, the following constrained approximation problem was studied: Minimize $F(\Phi)=\D(\Psi\|\Phi)=\int \log\left(\Psi/\Phi\right)\Psi$ among coercive spectra $\Phi\in L^\infty(\T)$ satisfying (\ref{moment}). Notice that minimization occurs with respect to the second argument. This permits to include the maximum entropy in this framework ($\Psi\equiv 1$) and to obtain a rational solution rather than one in the exponential class  when $\Psi$ is rational\footnote{Indeed, minimizing $\D(\Phi\|\Psi)$ with respect to $\Phi$ under (\ref{moment}) leads to extremal spectra of the form $\Phi_c=C\Psi\exp[G^\ast\Lambda G]$ which, due to the exponential factor, are non rational even when $\Psi$ is such.}. Further justification for this choice of the criterion may be found in \cite{GEORGIOU_LINDQUIST_KULLBACKLEIBLER}. In this case, for $\delta\Phi\in L^\infty$, $F'(\Phi;\delta\varphi)=-\langle\Phi^{-1}\Psi,\delta\varphi\rangle_{L^2}$. Since the constraint is as in (\ref{moment}), so is the space ${\mathcal V}^\perp$, see (\ref{ortcompl}). 
By Theorem \ref{hilbertspace}, we conclude that the optimal spectrum has the form obtained in  \cite{GEORGIOU_LINDQUIST_KULLBACKLEIBLER}
$$\Phi_c(e^{j\vartheta})=\frac{\Psi(e^{j\vartheta})}{G^\ast(e^{j\vartheta}) \Lambda_c G(e^{j\vartheta})},\quad \Lambda_c\in\Hermitian_n.
$$
The difficulties of extending this result to the multivariable case are illustrated in \cite[p.1062]{GEORGIOU_RELATIVEENTROPY}.

\section{Shannon entropy for finite measure spaces}\label{shannon}

The Shannon entropy underlying all the criteria so far considered will be here addressed directly via the first (rather than the second) part of equation (\ref{entropy}) and with a finite measure $\mu$ replacing Lebesgue measure. Let $(X, {\cal X},\mu)$ be a finite measure space and let $\varphi_i, i=1,\ldots,d$ be functions in ${\mathcal H}=L^2(X, {\cal X},\mu)$ and $\alpha\in\R^d$. Consider the problem of finding a nonnegative function $p$ in $L^\infty(X, {\cal X},\mu)$ maximizing the Shannon entropy
\begin{equation}\label{Sentropy}F(p)=H_\mu(p)=-\int_X\log [p(x)]p(x)d\mu
\end{equation}
under the constraints
\begin{eqnarray}\label{C1}&&\int_Xp(x)d\mu=1,\\&&\int_X\varphi_i(x)p(x)d\mu=\alpha_i,\quad i=1,\ldots,d.\label{C2}
\end{eqnarray}
Lemma \ref{chainr2}  can be readily adapted to this setting. Let $p_c\in L^\infty(X, {\cal X},\mu)$ be nonnegative and bounded away from zero $\mu$ a.e. Let $\delta p\in L^\infty(X, {\cal X},\mu)$. Then the directional derivative of the functional (\ref{Sentropy}) in direction $\delta p$ exists at $p_c$ and is given by
$$F'(p_c;\delta p)=\int_X \left[-1+\log p_c(x)\right]\delta p(x)d\mu=\langle -1+\log p_c,\delta p\rangle_{\mathcal H}.$$ 
Let us show that the fundamental geometric result Theorem \ref{hilbertspace} provides the form of the extremal solution also in this case. Suppose there exists $p_0\in L^\infty(X, {\cal X},\mu)$ a.e. everywhere positive satisfying (\ref{C1})-(\ref{C2}). Then $p\in L^\infty(X, {\cal X},\mu)$ also satisfies the constraints if it belongs to the affine space $p_0+{\cal V}$ where ${\cal V}$ is the subspace of functions $f\in L^\infty(X, {\cal X},\mu)$ such that
\begin{eqnarray}\label{S1}&&\int_Xf(x)d\mu=0,\\&&\int_X\varphi_i(x)f(x)d\mu=0,\quad i=1,\ldots,d.\label{S2}
\end{eqnarray}
Observe now that ${\cal V}^\perp$ is the subspace of functions of the form $\vartheta_0+\sum_{i=1}^d\vartheta_i\varphi_i(x)$. Observe also that for $p_c$ bounded and bounded away from zero as above, $\log p_c$ also belongs to $L^\infty(X, {\cal X},\mu)$ and, consequently, to $L^2(X, {\cal X},\mu)$.
By Theorem \ref{hilbertspace} we conclude that $(-1+\log p_c)\in {\mathcal V}^\perp$, it must namely be of the form  
\begin{equation}\label{exponential}p_c(x)=C\exp\left[\sum_{i=1}^d\vartheta_i\varphi_i(x)\right].
\end{equation}
If there exist constants $C$ and $\vartheta_i, i=1,\ldots,d$ such that $p_c$ belongs to $L^\infty(X, {\cal X},\mu)$, it is bounded away from zero $\mu$ a.e.  and it satisfies the constraints, then it is indeed optimal due to the concavity of the entropy. This is just the well-known fact that, if the maximizer exists, it belongs to the {\em exponential family}.
In the case when $d=1$ and  $\varphi_1=H$ the {\em Hamiltonian function}, we get a baby version of Gibbs variational principle, namely that the Gibbs distribution 
$$p_G(x)=C\exp\left[-\frac{H(x)}{kT}\right]$$ 
minimizes the {\em free energy}  $\langle H,p\rangle-kTF(p)$ where $F$ is as in (\ref{Sentropy}), $k$ is Boltzmann's constant and $T$ is absolute temperature. \cite{ellis}.

\section{Reciprocal processes identification with prior}\label{reciprocal}

In this section, we consider the problem of block-circulant covariance completion addressed in \cite{CFPP,CG} and we show that our result allows for a direct solution of this more general problem also in the case (not considered there) when a prior estimate is available.
The above mentioned block-circulant covariance completion  is equivalent to the computation of the parameters of a stationary, $m$-dimensional, {\em reciprocal process} of order $n$ defined on the discrete circle $\Z/N\Z$. A  process $y(t)$ defined on $\Z/N\Z$ taking values in $\R^m$ is reciprocal if it enjoys the following property. Take any two points $i,j\in\Z/N\Z$: They divide the discrete circle into two (discrete) arcs.  Then process $y(t)$ is reciprocal of order $1$ if $y(t)$ and $y(\tau)$ are conditionally independent given $y(i)$ and $y(j)$, for any $i,j$ and for any $t$ and $\tau$ belonging to different arcs. The process $y(t)$ is reciprocal of order $n$ if $y(t)$ and $y(\tau)$ are conditionally independent given $y(i),y(i+1), \dots y(i+n-1)$ and $y(j),y(j+1), \dots y(j+n-1)$, for any $i,j$ and for any $t$ and $\tau$ belonging to different arcs. 
Reciprocal processes  defined on (a  finite interval of) the integer line  can  be seen as a special class of discrete Markov random fields restricted to one dimension. Stationary reciprocal processes defined on $\Z/N\Z$ are potentially useful for describing  signals which naturally live in a finite region of the  time (or space) line such as texture images.

Let $\Sigma_i\in\R^{m\times m}$, $i=0,1,\dots, n$ be given.
In \cite{CFPP} the problem has been considered to compute the parameters of a stationary {\em reciprocal process} of order $n$ defined on the discrete circle $\Z/N\Z$ such that the first $n+1$ covariance lags of this process match the given  $\Sigma_i$, $i=0,1,\dots, n$.
For the importance and applications of this problem we refer to \cite{CFPP} and references therein. For a discussion of stationary reciprocal processes, we refer to 
\cite{Levy-F-98}.
In \cite{CFPP} is was shown that this problem is equivalent to
compute an extension $\Sigma_i\in\R^{m\times m}$, $i=n+1,n+2,\dots, N-1$ in such a way that the symmetric block-Toeplitz matrix $\Sigma \in\R^{Nm\times Nm}$ whose first block row is $[\Sigma_0\mid \Sigma_1\tp\mid \dots \Sigma_{N-1}\tp]$ maximizes 
\eq\label{recip-sp}
F(\Sigma):=\log[\det[\Sigma]]
\eeq
in the set ${\mathcal W}\cap {\mathcal S}$, where ${\mathcal S}$ is the cone of positive definite matrices and  ${\mathcal W}$ is the affine space of block-circulant symmetric matrices such that the north-west corner block of dimension
$m(n+1)\times m(n+1)$  is equal to the   symmetric block-Toeplitz matrix $\Sigma_{11}$ whose first block row is $[\Sigma_0\mid \Sigma_1\tp\mid \dots \Sigma_{n}\tp]$.
The form of solution to this problem may be easily computed by using Theorem \ref{theorem1}. In fact,
define
$$ 
U :=\bmat{ccccc}
0&I_m&0&\dots&0\\
0&0&I_m&\dots&0\\
\vdots&\vdots&
&\ddots&\vdots\\
       0&0&0&\dots&I_m\\
       I_m&0&0&\dots&0
       \emat\in\R^{Nm\times Nm},$$
       $$ E :=\bmat{ccccc} I_m&0&\ldots&0\\0&I_m&\ldots&0\\ 
					  0&0& \ddots&\vdots\\
				\vdots& &0&I_m\\
				0&0& \ldots&0  \emat\in \R^{ Nm\times (n+1)m}.
$$
where $I_m$ denotes the $m\times m$ identity matrix. Clearly, 
$U \tp U =U U \tp =I_{mN}$; i.e. $U $ is orthogonal. Note that a matrix $C$ with $N\times N$ blocks is   block-circulant  if and only if it commutes with $U $, namely if and only if it satisfies
\begin{equation}\label{charactcirc}
U \tp CU =C.
\end{equation}
The affine set ${\mathcal W}$ may be then characterized as 
\eq
{\mathcal W}=\{\Sigma=\Sigma\tp:\  E \tp\Sigma  E =\Sigma_{11},\ U \tp\Sigma U =\Sigma\}=A+{\mathcal V}
\eeq
with $A\in{\mathcal W}$ and
\eq
{\mathcal V}:=\{\Sigma=\Sigma\tp:\  E \tp\Sigma  E =0,\ U \tp\Sigma U =\Sigma\}.
\eeq
It is not difficult to check that
\eq
{\mathcal V}^{\perp}=\{\Delta= E \Lambda  E \tp+U \Theta U\tp-\Theta,\ 
\Lambda=\Lambda\tp\in\R^{(n+1)m\times (n+1)m},\ \Theta=\Theta\tp\in\R^{Nm\times Nm}\}.
\eeq
Hence the optimal solution, if it exists, has the form
\eq\label{optimsolr}
\Sigma_c=\left( E \Lambda  E \tp+U \Theta U\tp-\Theta\right)^{-1},\ 
\eeq
where $\Lambda=\Lambda\tp\in\R^{(n+1)m\times (n+1)m}$, and $\Theta=\Theta\tp\in\R^{Nm\times Nm}$
must be chosen in such a way that the constraints are satisfied.
This can be done through convex duality as discussed in \cite{CFPP}. The dual problem consists here in the unconstrained maximization of the concave function
$$ L(\Lambda,\Theta)=\Tr\log\left(E \Lambda  E \tp+U \Theta U\tp-\Theta\right)+ \Tr I-\Tr\left(\Lambda\Sigma_{11}\right).
$$
over a suitable set of multiplier pairs $(\Lambda,\Theta)$.
Once the optimal parameters $\Lambda$ and $ \Theta$ have been found, the optimal solution
(\ref{optimsolr})  has inverse $\Sigma_c^{-1}$ which is a block-circulant matrix whose first block-row has the form 
\begin{equation}\label{row}[M_0\mid M_1\mid\dots\mid M_n\mid0\mid0\mid\dots\mid 0\mid M_n\tp\mid
M_{n-1}\tp\mid\dots\mid M_1\tp],
\end{equation}
where the matrices $M_i$ are the sought for parameters of the stationary reciprocal process.

We now address the case when a prior information is available in terms of the parameters of a reciprocal process (possibly of higher order), or, equivalently of a prior positive definite covariance matrix $\Sigma_p\in\R^{Nm\times Nm}$.
In this case, instead of maximizing (\ref{recip-sp}) we minimize the divergence (see (\ref{divgauss}))
\eq\label{ind-gen-recip}
F(\Sigma):=\left[ \log\det
(\Sigma^{-1}\Sigma_p)+\tr(\Sigma_p^{-1}\Sigma)\right]
\eeq
under the same constraints.
By employing Theorem \ref{theoremmatrixprior}, we get the form of the optimal solution is
\eq\label{optimsolr-pr}
\Sigma_c=\left( E \Lambda  E \tp+U \Theta U\tp-\Theta+\Sigma_p^{-1}\right)^{-1},\ 
\eeq
where, again,  $\Lambda=\Lambda\tp\in\R^{(n+1)m\times (n+1)m}$, and $\Theta=\Theta\tp\in\R^{Nm\times Nm}$
must be chosen in such a way that the constraints are satisfied. As before, this can be done by solving a  dual problem for which existence can be proven along the lines of \cite{CFPP}. 
From (\ref{optimsolr-pr}) it follows that when $\Sigma_p$ is also the covariance matrix of a stationary reciprocal process of order $n$ or less, the optimal solution is also reciprocal of order $n$ and coincides with the optimal solution of the problem without prior.This result, analogous to what had been observed after Theorems \ref{theoremmatrixprior} and \ref{theorem4},  follows from (\ref{optimsolr-pr}) and the fact that there exists a unique block circulant covariance completion satisfying the linear constraints and having block  zeros in the first row as in (\ref{row}).
If instead,
$\Sigma_p$ is the covariance matrix of a stationary reciprocal process of order $n_1>n$ (requiring a larger memory),  
then the optimal solution is the covariance of a reciprocal process of order $n_1$ whose parameters may be read in the first block-row of $\Sigma_c^{-1}$.

\section{Extension to functionals defined on  a Banach space}\label{gen-ban}
In some applications, Theorem \ref{hilbertspace} does not suffice. 
For this reason, we mention the straightforward extension of our main result to functionals $F$ defined on a Banach space. 
 Let ${\mathcal X}$ be a Banach space and let  $F:{\mathcal X}\rightarrow \R$ be a functional. We say that $F$ is {\em G\^{a}teaux-differentiable} at $x_0$ in direction  $v$ if the limit
$$F'(x_0;v):=\lim_{\epsilon\rightarrow 0}\frac{F(x_0+\epsilon v)-F(x_0)}{\epsilon}
$$
exists. In this case, $F'(x_0;v)$ is called the {\em directional derivative} of $F$ at $x_0$ in direction $v$. We say that $F$ is {\em Fr\'echet-differentiable} at $x_0$  if there exists a bounded linear functional on ${\mathcal X}$ $DF_{x_0}$  such that
$$\lim_{\|x\|_{_{\mathcal X}}\rightarrow 0
}\frac{|F(x_0+x)-F(x_0)-DF_{x_0}(h)|}{\|x\|_{\mathcal X}}=0.
$$
The functional $DF_{x_0}$ is called the {\em Fr\'echet differential} of $F$ at $x_0$. Again, if $F$ is Fr\'echet differentiable at $x_0$, then $DF_{x_0}$ is unique and, for any $x\in{\mathcal X}$, $F$ is G\^{a}teaux differentiable at $x_0$ in direction $v$ and it holds
\begin{equation}\label{FGrelBanach}
F'(x_0;v)=DF_{x_0}(v).
\end{equation}

\begin{theorem}\label{banachspace}Let ${\mathcal X}$ be a Banach space, let ${\mathcal V}\subseteq {\mathcal X}$ be a subspace, let $x\in{\mathcal X}$ and consider the corresponding {\em coset} ${\mathcal W}:= x+{\mathcal V}$. Assume that the functional $F$ is Fr\'echet-differentiable  at $w_c\in{\mathcal W}$.
Then $w_c$ is a critical point of $F$ over ${\mathcal W}$ if and only if $ DF_{w_c}$ belongs to the {\em annihilator} of ${\mathcal V}$.
\end{theorem}
\proof
Observe that $F'(w_c;v)=0$ for all $v\in{\mathcal V}$ if and only if $DF_{w_c}(v)=0, \forall v\in{\mathcal V}$.
\qed
When $F$ is not Fr\'echet-differentiable  at $w_c$ but merely G\^{a}teaux differentiable in directions varying in a subspace, a generalization such as in Theorem \ref{hilbertspace} can be established. Nevertheless, we like to give here an even more general result which can be effectively applied when the solution lies on the boundary of the feasible set. Indeed, all the maximum entropy applications so far considered in this paper feature a solution which is an {\em interior point} of the admissible set \cite{BYRNES_LINDQUIST_INTERIOR}. Let us begin by recalling the fundamental result of convex optimization. Let $K$ be a convex subset of the vector space $X$, let $F:K\rightarrow \R$ be  concave and let $x_0\in K$. Then, the one-sided directional derivative or hemidifferential of $F$ at $x_0$ in direction $x-x_0$ 
$$F'_+(x_0;x-x_0):=\lim_{\epsilon\searrow 0}\frac{F(x_0+\epsilon(x-x_0))-F(x_0)}{\epsilon}
$$
exists for every $x\in K$ (this is a consequence of the monotonicity of the difference quotients) \cite[p.66]{K}.
\begin{theorem}Let $K$ be a convex subset of the vector space $X$ and  let $F:K\rightarrow \R$ be  concave. Then, $x_0\in K$ is a maximum point for $F$ over $K$ if and only if it holds
\begin{equation}\label{suffopt}
F'_+(x_0;x-x_0)\le 0, \quad \forall x\in K.
\end{equation}
\end{theorem}
As a corollary, we get the following {\em sufficient} condition for optimality.
\begin{corollary}\label{corsuff}Let ${\mathcal X}$ be a Banach space, let ${\mathcal V}\subseteq {\mathcal X}$ be a subspace, let $x\in{\mathcal X}$ and consider the corresponding {\em coset} ${\mathcal W}:= x+{\mathcal V}$. Let $K$ be a convex subset of ${\mathcal W}$ and let $F:K\rightarrow \R$ be concave. Assume that, for any $x\in K$, the hemidifferential of $F$ at $x_c\in K$ in direction $x-x_c$ is given through the linear continuous functional $D_K F(x_c)\in X^\ast$ as
\begin{equation}\label{reprhemi}
F'_+(x_c;x-x_c)=\langle D_K F(x_c),x-x_c\rangle,
\end{equation}
where $\langle \cdot,\cdot\rangle$ is the duality pairing between $X^\ast$ and $X$. Then, if $D_K F(x_c)$ belongs to the annihilator of ${\mathcal V}$, $x_c$ is a maximum point of $F$ over $K$.
\end{corollary}
\proof
Observe that $(x-x_c)\in {\mathcal V}, \forall x\in K$. If $D_K F(x_c)$ belongs to the annihilator of ${\mathcal V}$, by (\ref{reprhemi}), we get $F'_+(x_0;x-x_0)=0$, for all $x\in K$. By (\ref{suffopt}), $x_c$ is optimal.
\qed
This result permits to establish optimality of the solutions computed in Sections \ref{matricial functions} and \ref{prior}  in a  larger class of spectra.

As a simple application of Corollary \ref{corsuff}, we now show that the Gaussian has maximum entropy
among all probability densities with given mean and variance. Let $X=L^1(\R)$ and consider the affine space ${\cal W}$ of $L^1$-functions $f$ satisfying the constraints
\begin{equation}\label{3constr}\int_\R f(x)dx=1,\quad \int_\R xf(x)dx=0,\quad \int_\R x^2f(x)dx=\sigma^2.
\end{equation}
Observe, as in Section \ref{shannon}, that the corresponding subspace ${\mathcal V}$  is given by  $L^1$-functions $f$ satisfying
\begin{equation}\label{3conV}\int_\R f(x)dx=0,\quad \int_\R xf(x)dx=0,\quad \int_\R x^2f(x)dx=0.
\end{equation}

Let $K$ be the convex subset of ${\cal W}$ obtained by intersecting ${\mathcal W}$ with the cone of nonnegative functions $p$. Let us take  as criterion on $K$ the concave functional given by the Shannon entropy
$$F(p)=H(p)=-\int_X\log [p(x)]p(x)dx.
$$
The hemidifferential of $H(p)$ at $p_c\in K$ in direction $p-p_c, p\in K$, has the form
\begin{eqnarray}\label{formhemidiff}H'(p_c;p-p_c)=\int_\R \left[-1+\log p_c(x)\right](p(x)-p_c(x))dx\\=\langle -1+\log p_c,p-p_c\rangle=\langle D_K H(p_c),p-p_c\rangle.
\end{eqnarray}
If $p_c$ has the form
$$p_c(x)=C\exp\left[\vartheta_1 x+\vartheta_2 x^2\right],
$$
then,  in view of (\ref{3conV}) and (\ref{formhemidiff}), $D_K H(p_c)$ belongs to the annihilator of ${\mathcal V}$. By Corollary \ref{corsuff}, such a $p_c$ is optimal provided it belongs to $K$. For $C=(2\pi)^{-1/2}$, $\vartheta_1=0$ and $\vartheta_2=-\frac{1}{2\sigma^2}$, we get that $p_c$ is nonnegative and satisfies the constraints (\ref{3constr}), i.e. it belongs to $K$.  Thus, the Gaussian density 
$$p_c(x)=(2\pi)^{-1/2}\exp \left[-\frac{1}{2}\frac{x^2}{\sigma^2}\right]$$ 
has maximum entropy among densities with given mean and variance.

\section{Closing comments}
In this paper, we have established a simple orthogonality condition that allows to derive the form of the optimal solution in a plethora of maximum entropy problems.
We feel that this geometric condition affords a considerable conceptual simplification allowing to cast least-squares and maximum entropy problems in the same framework (admittedly, not as deep as the one provided in \cite{csiszar2}). It can, moreover, be readily generalized to abstract situations and to problems with nonlinear constraints. Further study is needed to see whether  this approach may be suitably adapted to the abstract setting of Subsection \ref{schroedinger}. A suitable mixture of the geometry we have seen in Burg's and in Dempster's problems in Subsections \ref{burgME} and \ref{dempstermaxent} might provide the key to understanding AR and ARMA Identification of Graphical Models, a topic which has recently received considerable attention, see e.g. \cite{lauritzen,dahlhaus,HL,vandenberghe,vandenberghe2,avventi}. Finally, we should never forget the motto over the entrance to Plato's Academy:  ``$A\gamma\epsilon\omega\mu\acute{\epsilon}\tau\rho\eta\tau o\varsigma\; \mu\eta\delta\epsilon\grave{\iota}\varsigma \;\epsilon\iota\sigma\acute{\eta}\tau\omega$", namely ``Let no one untrained in geometry enter."

\section*{Acknowledgments}
The authors wish to thank two anonymous reviewers for a careful reading and for providing several constructive suggestions.  In particular, we are thankful to one reviewer for suggesting to employ a generalization without  Fr\'echet differentiability of  the main result, for encouraging us to access Boltzmann's original work \cite{BOL} and for providing an endless string of technical and expository suggestions that led to a considerable improvement of the paper.

\end{document}